\documentclass[a4paper,10pt]{amsart}
\usepackage[arrow,matrix]{xy}
\usepackage{amsmath,amssymb,amscd,bbm,amsthm,mathrsfs}
\usepackage[colorlinks=true,citecolor=blue,linkcolor=blue]{hyperref}

\theoremstyle{plain}
\newtheorem{thm}{Theorem}[section]
\newtheorem{cor}[thm]{Corollary}
\newtheorem{lem}[thm]{Lemma}
\newtheorem{prop}[thm]{Proposition}

\newtheorem{exa}[thm]{Example}

\newtheorem{rem}[thm]{Remark}

   \def\op{\oplus} \def\ot{\otimes}
\def\Hom{\operatorname {Hom}}

\def\Ext{\operatorname {Ext}} 
 \def\k{\mathbbm{k}}

\begin{document}
\title[Skew polynomial algebras]{\bf Skew polynomial algebras with coefficients in Koszul Artin-Schelter regular algebras}

\author{Ji-Wei He, Fred Van Oystaeyen and Yinhuo Zhang}
\address{J.-W. He\newline \indent Department of Mathematics, Shaoxing College of Arts and Sciences, Shaoxing Zhejiang 312000,
China\newline \indent Department of Mathematics and Computer
Science, University of Antwerp, Middelheimlaan 1, B-2020 Antwerp,
Belgium} \email{jwhe@usx.edu.cn}
\address{F. Van Oystaeyen\newline\indent Department of Mathematics and Computer
Science, University of Antwerp, Middelheimlaan 1, B-2020 Antwerp,
Belgium} \email{fred.vanoystaeyen@ua.ac.be}
\address{Y. Zhang\newline
\indent Department WNI, University of Hasselt, Universitaire Campus,
3590 Diepenbeek, Belgium} \email{yinhuo.zhang@uhasselt.be}

\date{}
\begin{abstract}
Let $A$ be a Koszul Artin-Schelter regular algebra with Nakayama automorphism $\xi$. We show that the Yoneda Ext-algebra of the skew polynomial algebra $A[z;\xi]$ is a trivial extension of a Frobenius algebra. Then we prove that $A[z;\xi]$ is Calabi-Yau; and hence each Koszul Artin Schelter regular algebra is a subalgebra of a Koszul Calabi-Yau algebra. A superpotential $\widehat{w}$ is also constructed so that the Calabi-Yau algebra $A[z;\xi]$ is isomorphic to the derivation quotient of $\widehat{w}$. The Calabi-Yau property of a skew polynomial algebra with coefficients in a PBW-deformation of a Koszul Artin-Schelter regular algebra is also discussed.
\end{abstract}


\maketitle

\section*{Introduction}

Let $\k$ be a field of characteristic zero. Let $A=\op_{k\in \Bbb{Z}}A_k$ be a ($\Bbb{Z}$-)graded ($\k$-)algebra, and $M=\op_{k\in \Bbb{Z}}M_k$ be a graded left $A$-module. The $n^{th}$ shift of $M$ is the graded $A$-module $M(n)$ whose $k^{th}$ component is: $M(n)_k=M_{n+k}$. If $M$ is a graded $A$-bimodule and $\sigma,\varphi$ are graded automorphisms of $A$, then ${}_\sigma M_\varphi$ is the graded $A$-bimodule whose left $A$-action is twisted by $\sigma$ and right $A$-action is twisted by $\varphi$. A graded algebra $A$ is called {\it Calabi-Yau} of dimension $d$, if (cf. \cite{Gin}):

(i) $A$ is
homologically smooth; that is, $A$ has a bounded resolution of
finitely generated graded projective $A$-bimodules;

(ii) $\Ext^i_{A^e}(A,A\ot A)=0$ if $i\neq d$ and $\Ext_{A^e}^d(A,A\ot A)\cong A(l)$ for some integer $l$ as $A$-bimodules, where $A^e=A\ot A^{op}$ is the enveloping algebra of $A$.

Calabi-Yau algebras are strongly related to Artin-Schelter (AS, for short) regular algebras. Recall that a connected graded algebra $A=\k\op A_1\op A_2\op\cdots$ is  called an {\it AS-regular} algebra if (i) $A$ has finite global dimension $d$; (ii) $\underline{\Ext}_A^i({}_A\k,A)=0$ if $i\neq d$, and $\dim\underline{\Ext}_A^d({}_A\k,A)=1$. Here $\underline{\Ext}$ is the derived functor of graded $\underline{\Hom}$ (cf. \cite{Smi}). If a graded Calabi-Yau algebra is also connected, then it must be AS-regular \cite{BT}. If an AS-regular algebra $A$ of global dimension $d$ is Noetherian or Koszul, then $A$ differs from a Calabi-Yau algebra by an automorphism; more precisely, we have $\Ext^i_{A^e}(A,A\ot A)=0$ for $i\neq d$, and $\Ext^d_{A^e}(A,A\ot A)\cong A_\xi(l)$ for some graded automorphism $\xi$ of $A$ \cite{VdB,BM}. We call the automorphism $\xi$ the {\it Nakayama automorphism} of $A$.

Berger and Pichereau recently constructed in \cite{BP} an interesting class of Calabi-Yau algebras of dimension 3, which are related to deformations of Poisson algebras. Given an AS-regular algebra $A$ of global dimension 2 (which must be Koszul), Dubois-Violette showed in \cite{DV} (also see \cite{Z})  that $A$ is defined by an invertible matrix $M$, that is, $A\cong T(V)/(f)$ with $V$ a finite dimensional vector space with a basis $\{x_1,\dots,x_n\}$ and $f=(x_1,\dots,x_n)M(x_1,\dots,x_n)^t$. Here the matrix multiplications should be regarded in $T(V)$. Given such an AS-regular algebra, Berger and Pichereau constructed a graded algebra $B(f)$ which is generated by $V\op\k z$, and whose generating relations are cyclic partial derivations of $w=fz$. They proved that $B(f)$ is a Calabi-Yau algebra of dimension 3, and gave the classification, up to isomorphisms, of the obtained Calabi-Yau algebras. They also showed that $B(f)$ is isomorphic to the skew polynomial algebra $A[z;\xi]$ for some automorphism $\xi$ of $A$. We find that $\xi$ is exactly the Nakayama automorphism of $A$, and the Calabi-Yau property of $A[z;\xi]$ holds for general Koszul AS-regular algebra $A$ by inspecting the Yoneda-Ext algebra of $A[z;\xi]$. The main results of this paper are the following (cf. Theorem \ref{thm1}).

\begin{thm} Let $A$ be a Koszul AS-regular algebra of global dimension $d$, and $\xi$ the Nakayama automorphism of $A$. Then
the skew polynomial algebra $B=A[z;\xi]$ is a Calabi-Yau algebra of dimension $d+1$.
\end{thm}
Since an AS-regular algebra of global dimension 2 is always Koszul, our main results provide a new proof of \cite[Theorem 2.10]{BP} and our proof is totally different from that in \cite{BP}.

It was shown in \cite{DV} and \cite{BSW} that a Koszul AS-regular algebra $A$ is determined by a twisted superpotential $w$. We show that the twisted superpotential $w$ can be symmetrized into a superpotential $\widehat{w}$ by introducing a new indeterminate, so that the skew polynomial algebra $A[z;\xi]$ is isomorphic to the derivation quotient algebra obtained from the superpotential $\widehat{w}$ (see Theorem \ref{thm2}).

Let $A$ be a Koszul AS-regular algebra, and $\xi$ the Nakayama automorphism of $A$. A PBW-deformation of $A$ is a filtered algebra $U$ such that its associated graded algebra is isomorphic to $A$. For a PBW-deformation $U$ of $A$, $U$ has a filtration-preserving automorphism $\zeta$ such that $gr(\zeta)=\xi$, still called a Nakayama automorphism (in this case, $\zeta$ is not unique, see more details in Section \ref{pbw}). It is natural to ask whether $U[z;\zeta]$ is a Calabi-Yau algebra. Recall that a nongraded algebra $U$ is {\it Calabi-Yau} of dimension $d$ if (i) $U$ is
homologically smooth; (ii) $\Ext^i_{U^e}(U,U\ot U)=0$ if $i\neq d$ and $\Ext_{U^e}^d(U,U\ot U)\cong U$ as $U$-bimodules.

Since a Nakayama automorphism $\zeta$ respects the filtration of $U$, we see that $U[z;\zeta]$ is in fact a PBW-deformation of $A[z;\xi]$, which is a Koszul algebra. Then we can use the techniques developed in \cite{PP,Po} for PBW-deformations of Koszul algebras to discuss the Calabi-Yau property of $U[z;\zeta]$.

Now let $A$ be a Koszul Calabi-Yau algebra, and $U$ be a PBW-deformation of $A$. We may choose a specific Nakayama automorphism $\zeta$ of $U$ (see Proposition \ref{prop3}) so that we have (cf. Theorem \ref{thm4}):

\begin{thm} $U[z;\zeta]$ is Calabi-Yau.
\end{thm}

In the theorem above, if $A$ is only an AS-regular algebra, then the result may fail. Counterexamples may be found in the case where $A$ is AS-regular of global dimension 2. At the end of the paper, we provide a necessary and sufficient condition for $U[z;\zeta]$ to be Calabi-Yau with $U$ a PBW-deformation of an AS-regular algebra of global dimension 2 (cf. Theorem \ref{thm5}).

\section{Trivial extensions}\label{sec1}

Let $E=\k\op E_1\op E_2\op\cdots$ be a connected graded algebra, and $M$ a graded $E$-bimodule. Recall that the trivial extension of $E$ by $M$ is the graded algebra $\Gamma(E,M)=E\op M$ with the product $(x_1,m_1)*(x_2,m_2)=(x_1x_2,x_1\cdot m_2+m_1\cdot x_2)$ for $x_i\in E$ and $m_i\in M$. If
$M_i=0$ for all $i\leq0$, then $\Gamma(E,M)$ is a connected graded algebra with the $i^{th}$ component $\Gamma(E,M)_i=E_i\op M_i$.

We focus on trivial extensions of finite dimensional algebras. Let $E$ be a finite dimensional connected graded algebra. We say that $E$ is of {\it length} $d$ if $E_d\neq0$ and $E_i=0$ for all $i>d$. Let $E$ be a connected finite dimensional algebra of length $d$, and $\sigma$ a graded automorphism of $E$. Let $E^*$ be the dual vector space of $E$. Then $E^*$ is a graded $E$-bimodule with the induced $E$-action. Let $E^*_\sigma$ be the graded $E$-bimodule obtained from $E^*$ with the right $E$-action twisted by $\sigma$. Given an integer $n>d$, consider the trivial extension of $E$ by the bimodule $E^*_{\sigma}(-n)$: $\Gamma(E,E^*_\sigma(-n))= E\op E^*_\sigma(-n)$. For simplicity, we write $\Gamma(E,\sigma,n)$ for $\Gamma(E,E^*_\sigma(-n))$. Now the product of $\Gamma(E,\sigma,n)$ is defined by: $(x_1,f_1)*(x_2,f_2)=(x_1x_2,x_1\cdot f_2+f_1\cdot\sigma(x_2))$ for $x_1,x_2\in E$ and $f_1,f_2\in E^*$.

For later discussions, we introduce first some terminology. A connected graded algebra $E$ of length $d$ is called a graded Frobenius algebra if there is an isomorphism of graded left $A$-modules $\Theta:E\cong E^*(-d)$, or equivalently, there is a nondegenerated graded bilinear form $\langle\ ,\ \rangle:E\times E\to \k(-d)$ such that $\langle x,yz\rangle=\langle xy,z\rangle$ for all $x,y,z\in E$. If $E$ is graded Frobenius, then there is a unique graded automorphism $\varphi$ of $E$ such that $\Theta:E_\varphi\to E^*(-d)$ is an isomorphism of $A$-bimodules. The isomorphism $\varphi$ is called the {\it Nakayama automorphism} of $E$. For the bilinear form, we have $\langle x,y\rangle=\langle y,\varphi(x)\rangle$. A graded Frobenius algebra $E$ of length $d$ is said to be {\it graded symmetric}, if $\langle x,y\rangle=(-1)^{i(d-i)}\langle y,x\rangle$ for all $x\in E_i$ and $y\in E_{d-i}$. In this case, $\varphi=\epsilon^{d-1}$, where $\epsilon:E\to E$ is defined by $\epsilon(x)=(-1)^ix$ for $x\in E_i$.

\begin{prop}\label{trivial-ext} Let $E$ be a connected graded algebra of length $d$, and $\sigma$ a graded automorphism of $E$. Then the trivial extension $\Gamma(E,\sigma,n)$ $(n>d)$ is a graded Frobenius algebra of length $n$, and the Nakayama automorphism $\varphi$ of $\Gamma(E,\sigma,n)$ is given by $\varphi(x,f)=(\sigma^{-1}(x),f\circ\sigma)$ for all $x\in E$ and $f\in E^*$.
\end{prop}
\proof Define a bilinear form $\langle\ ,\ \rangle:\Gamma(E,\sigma,n)\times\Gamma(E,\sigma,n)\to \k(-n)$ by $\langle(x_1,f_1),(x_2,f_2)\rangle=f_2(x_1)+f_1(\sigma(x_2))$ for $x_1,x_2\in E$ and $f_1,f_2\in E^*$. A straightforward verification shows that $\langle(x_1,f_1),(x_2,f_2)*(x_3,f_3)\rangle=\langle(x_1,f_1)*(x_2,f_2),(x_3,f_3)\rangle$. Obviously, the bilinear form $\langle\ ,\ \rangle$ is nondegenerated. Hence $\Gamma(E)$ is a Frobenius algebra. Moreover, $\langle(x_1,f_1),(x_2,f_2)\rangle=f_2(x_1)+f_1(\sigma(x_2))=f_2(\sigma\circ\sigma^{-1}(x_1))+f_1\circ\sigma(x_2)
=\langle(x_2,f_2),(\sigma^{-1}(x_1),f_1\circ\sigma)\rangle$. Hence the Nakayama automorphism $\varphi$ of $\Gamma(E,\sigma,n)$ is defined by $\varphi(x,f)=(\sigma^{-1}(x),f\circ\sigma)$ for all $x\in E$ and $f\in E^*$. \qed

\begin{rem}\label{rem1} {\rm In Proposition \ref{trivial-ext}, if we choose $\sigma=\epsilon^{n-1}$, then the Nakayama automorphism $\varphi$ is given as follows: for $x\in E_i$ and $f\in E^*_i$, $\varphi(x,f)=((-1)^{(n-1)i}x,f\circ\epsilon^{(n-1)i})=((-1)^{(n-1)i}x,(-1)^{(n-1)i}f)$. Hence $\varphi=\epsilon^{n-1}$. Therefore, $\Gamma(E,\epsilon^{n-1},n)$ is a graded symmetric algebra.}
\end{rem}

Recall that a (cochain) differential graded algebra (dga, for short) $(E,\delta_E)$ is a graded algebra $E=\op_{n\in \Bbb Z}E_n$ together with a derivation $\delta_E$ of degree 1 such that $\delta^2_E=0$. A left differential graded module ${}_EM$ is a left graded $E$-module together with a differential $\delta_M$ such that $\delta_M(xm)=\delta_E(x)m+(-1)^{|x|}x\delta_M(m)$ for all homogeneous elements $x\in E$ and $m\in M$, where $|x|$ denotes the degree of $x$. Similarly, one has right differential graded modules and differential graded bimodules.

A {\it curved differential graded algebra} (cdga, for short) (cf. \cite{PP, Po}) is a triple $(E,\delta_E,\theta_E)$, where $E$ is a graded algebra, $\delta_E$ is a derivation of degree 1 and $\theta_E$ is a special element in $E_2$, such that $\delta_{E}(\theta_E)=0$ and $\delta_E^2(x)=\theta_E x-x\theta_E$ for all homogeneous elements $x\in E$. The element $\theta_E$ is usually called the {\it curvature element} of $E$. Let $(E',\delta_{E'},\theta_{E'})$ be another cdga. A {\it cdga morphism} $f:E\to E'$ is a graded algebra morphism such that $f(\theta_E)=\theta_{E'}$ and $f\delta_{E}=\delta_{E'}f$ (warning: our definition of cdga morphism given here is more restricted than that in \cite{PP}). A {\it cdg $E$-bimodule} is a graded $E$-bimodule $M$ endowed with a differential $\delta_M$ which is compatible with the differential $\delta_E$ of $E$ and satisfies the condition $\delta^2_M(m)=\theta_E m-m\theta_E$. Note that if the curvature element is zero, then a cdga is just a usual dga, and a cdg bimodule is a usual dg bimodule.

Let $(E,\delta_E,\theta_E)$ be a cdga, and let $M$ be a cdg $E$-bimodule. The trivial extension of $E$ by $M$ is the cdga $(\Gamma_{cdg}(E,M),\delta_{\Gamma_{cdg}},\theta_{\Gamma_{cdg}})$ defined as follows: as a graded algebra $\Gamma_{cdg}(E,M)$ is just the trivial extension $\Gamma(E^\sharp,M^\sharp)$, where $E^\sharp$ is the underlying graded algebra by forgetting the derivation $\delta_E$ of $E$, and $M^\sharp$ is the underlying graded bimodule of $M$; the derivation $\delta_{\Gamma_{cdg}}$ is defined by $$\delta_{\Gamma_{cdg}}(x,m)=(\delta_E(x),\delta_M(m))$$ for all $x\in E$ and $m\in M$; and the curvature element $\theta_{\Gamma_{cdg}}=(\theta_E,0)$.

Let $M(n)$ be the $n^{th}$ shift of the graded $E$-bimodule $M$. Note that the differential $\delta_{M(n)}$ and the $E$-actions of $M(n)$ should be changed slightly so that $M(n)$ is also a cdg $E$-bimodule: For a homogeneous element $m\in M$, we denote by $m(n)$ the corresponding element in $M(n)$. Then $\delta_{M(n)}(m(n))=(-1)^n\delta_M(m)(n)$. Let $x\in E$ be a homogeneous element. The left $E$-action on $M(n)$ is defined by $x\diamond (m(n))=(-1)^{n|x|}(x\cdot m)(n)$ and the right $E$-action is defined by $(m(n))\diamond x=(m\cdot x)(n)$, where $x\cdot m$ and $m\cdot x$ are $E$-actions on $M$.

Let $M^\vee=\op_{n\in \Bbb Z}M_n^*$ be the graded dual of $M$. Then $M^\vee$ is a cdg $E$-bimodule with the differential $\delta_{M^\vee}$ and $E$-actions defined as follows: for homogeneous elements $f\in M^\vee$, $m\in M$ and $x\in E$, we have \begin{equation}\label{eq6}
    \delta_{M^\vee}(f)=(-1)^{|f|+1}f\circ \delta_M,
\end{equation}
\begin{equation}\label{eq7}
    (x\rightharpoonup f)(m)=(-1)^{|x|(|f|+|m|)}f(m\cdot x)\text{ and }(f\leftharpoonup x)(m)=f(x\cdot m).
\end{equation}

Now let $E=\k\op E_1\op\cdots \op E_d$ $(E_d\neq0)$ be a finite dimensional cdga with differential $\delta_E$. Then $E^*=E^\vee$ is a cdg $E$-bimodule. Hence the trivial extension $\Gamma_{cdg}(E,E^*(-d-1))$ is a cdga.

Let $M$ be a cdg $E$-bimodule. Note that $(M(n))^\sharp$ is different from $M^\sharp(n)$ as graded $E^\sharp$-bimodules. We remark that the graded algebra $\Gamma_{cdg}(E,E^*(-d-1))^\sharp$ is different from the trivial extension $\Gamma(E^\sharp,(E^*)^\sharp(-d-1))$ of the graded algebra $E^\sharp$. However, we have the following result.

\begin{lem}\label{lem4} As a graded algebra, we have $$\Gamma_{cdg}(E,E^*(-d-1))^\sharp= \Gamma(E^\sharp,{}_{\epsilon^d}(E^*)^\sharp(-d-1)).$$
\end{lem}
\proof As before, we denote by $x\cdot f$ and $f\cdot x$ for $x\in E$ and $f\in E^*$ the $E^\sharp$-actions on the graded dual $(E^*)^\sharp$. Note that in the definition (\ref{eq7}) of the $E$-actions on $E^*$, $f(x\cdot m)$ is zero unless $|x|+|m|=-|f|$. Hence we have $(x\rightharpoonup f)=(-1)^{|x|}x\cdot f$.

It suffices to verify the multiplication of the cdga $\Gamma_{cdg}(E,E^*(-d-1))$ is equal to the multiplication of $\Gamma(E^\sharp,{}_{\epsilon^d}(E^*)^\sharp(-d-1))$. Indeed, for homogeneous elements $x,y\in E$ and $f,g\in E^*$, we have
$$\begin{array}{ccl}
    (x,f(d+1))*(y,g(-d-1))&=&\left(xy,x\diamond(g(-d-1))+(f(-d-1))\diamond y\right)\\
    &=&\left(xy,(-1)^{|x|(d+1)}(x\rightharpoonup g)(-d-1)+(f\leftharpoonup y)(-d-1)\right)\\
    &=&\left(xy,(-1)^{|x|d}(x\cdot g)(-d-1)+(f\cdot y)(-d-1)\right)\\
    &=& \left(xy,(\epsilon^{d}(x)\cdot g)(-d-1)+(f\cdot y)(-d-1)\right)\\
  \end{array}
$$ Now it is easy to see that the last item in the identities above is exactly the multiplication of elements $(x,f(d+1))$ and $(y,g(-d-1))$ in the trivial extension $\Gamma(E^\sharp,{}_{\epsilon^d}(E^*)^\sharp(-d-1))$. \qed

\section{Yoneda algebras}\label{sec}

In this section, we will compute the Yoneda products of a skew polynomial algebra with coefficients in a Koszul algebra. We first recall the definition of a Koszul algebra. Let $V$ be a finite dimensional vector space. A {\it quadratic} algebra $A$ is a connected graded algebra of form $A=T(V)/(R)$, where $R\subseteq V\ot V$ and $(R)$ is the two-sided ideal of $T(V)$ generated by $R$. The {\it quadratic dual} $A^!$ of a quadratic algebra $A$ is defined to be $A^!=T(V^*)/(R^\perp)$, where $R^\perp\subseteq V^*\ot V^*$ is the orthogonal complement of $R$. One easily sees that $(A^!)^!= A$. Let $\phi:A\to A$ be an automorphism of the quadratic algebra $A$. The restriction of $\phi$ to $A_1=V$ induces a bijective linear map $f:V^*\to V^*$. Since $A$ is quadratic, we see that $f$ defines an automorphism $\phi^!$ of the quadratic dual algebra $A^!$. We call $\phi^!$ the {\it automorphism of $A^!$ dual to} $\phi$. Since $(A^!)^!= A$, we have $(\phi^!)^!=\phi$.

A quadratic algebra $A$ is called a {\it Koszul algebra} \cite{Pr} if the trivial graded module ${}_A\k$ admits a graded projective resolution: $$\cdots\longrightarrow P^{-n}\longrightarrow\cdots\longrightarrow P^{-1}\longrightarrow P^0\longrightarrow{}_A\k\longrightarrow0,$$ such that the graded module $P^{-n}$ is generated in degree $n$ for all $n\ge0$. Recall that if $A$ is Koszul, then the Yoneda algebra $E(A):=\op_{i\ge0}\Ext_A^i({}_A\k,{}_A\k)\cong A^!$. Moreover, $A$ is Koszul if and only if $A^!$ is Koszul \cite{Smi}. We refer to \cite{Pr} and \cite{Smi} for further properties of Koszul algebras.

Let $A=T(V)/(R)$ be a Koszul algebra, and let $C_0=\k$, $C_{-1}=V$, $C_{-2}=R$ and $C_{-n}=\bigcap_{i=0}^{n-2}V^{\ot i}\ot R\ot V^{n-i-2}$ for $n\ge3$. Then $C=\bigoplus_{n\ge0}C_{-n}$ is a graded subcoalgebra of the tensor coalgebra $T(V)$. Moreover, as graded algebras, $E(A)\cong C^\vee=\bigoplus_{n\ge0}C^*_{-n}\cong A^!$.

Consider the graded minimal projective resolution of the trivial module ${}_A\k$:
\begin{equation}\label{minres}
    \cdots \longrightarrow A\ot C_{-n}\overset{\partial^{-n}}\longrightarrow\cdots
\overset{\partial^{-2}}\longrightarrow A\ot C_{-1}\overset{\partial^{-1}}\longrightarrow A\longrightarrow{}_A\k\longrightarrow0,
\end{equation} where the differential is given on pure tensors by:
$$   \partial^{-n}(a\ot x_1\ot\cdots\ot x_n) =  ax_1\ot x_2\ot\cdots\ot x_n, $$ for all $a\in A$ and $x_1,\dots,x_n\in V$.

Let $\sigma$ be a graded automorphism of $A$. Since $A$ is Koszul, $\sigma$ induces an automorphism (also denoted by $\sigma$) of $C$ in the obvious way. Let $\sigma^\vee$ be the automorphism of graded algebra $C^\vee$ induced by $\sigma$. Since $A^!\cong C^\vee$, we see that $\sigma^\vee=\sigma^!$. Let $B=A[z;\sigma]$ be the algebra of skew polynomials with coefficients in $A$. We assume that $z$ is of degree 1. Then it is well known that $B$ is also a Koszul algebra (cf. \cite{ST}, for example). The elements of $B$ are of the sums of the elements of the form $az^i$ with $i\ge0$ and $a\in A$, moreover $za=\sigma(a)z$. We want to construct a minimal projective resolution of the trivial module ${}_B\k$. The following construction is standard (cf. \cite{GS} or \cite{P}).

Clearly, $B$ is free both as a left $A$-module or as a right $A$-module. Applying the exact functor $B\ot_A-$ to the projective resolution (\ref{minres}) of $\k$, we obtain the following complex:
\begin{equation}\label{Btensor}
    \cdots \longrightarrow B\ot C_{-n}\longrightarrow\cdots
\longrightarrow B\ot C_{-1}\longrightarrow B\longrightarrow0.
\end{equation} The complex is exact except at the final position. The cohomology at the final position is $B/BA_+$. By abusing the notation, we also denote the differential of the complex (\ref{Btensor}) by $\partial$.

For each $n\ge1$, we define a homomorphism of left $B$-modules $$f^{-n}:B\ot C_{-n}\longrightarrow B\ot C_{-n}$$ by $$f^{-n}(1\ot x_1\ot\cdots\ot x_n)=z\ot \sigma^{-1}(x_1)\ot\cdots\sigma^{-1}(x_n)$$ for all $x_1,\dots,x_n\in V$. In addition, define a left $B$-module homomorphism $f^0:B\to B$ by $f^0(1)=z$. It is easy to check that these $f^{-n}$ are compatible with the differential of the complex (\ref{Btensor}). Hence $f=\prod_{n\ge0} f^{-n}$ is a morphism of complexes. The mapping cone of $f$ reads as follows:
{\small\begin{equation}\label{Bres}
 \cdots \longrightarrow B\ot C_{-n}\op B\ot C_{-n+1}\overset{\delta^{-n}}\longrightarrow\cdots
\longrightarrow B\ot C_{-2}\op B\ot C_{-1}\overset{\delta^{-2}}\longrightarrow B\ot C_{-1}\op B\overset{\delta^{-1}}\longrightarrow B\longrightarrow0,
\end{equation}} where the differential is given by:
$\delta^{-n}=\left(
   \begin{array}{cc}
     \partial^{-n} & f^{-n+1} \\
     0 & -\partial^{-n+1} \\
   \end{array}
 \right)
$ for $n\ge2$, and $\delta^{-1}=(\partial^{-1},f_0)$ for $n=1$.

A straightforward verification shows that the complex (\ref{Bres}) is exact except at the final position, and the cohomology at the final position is $\k$. Hence we have the following lemma.

\begin{lem}\label{lem-res} The complex (\ref{Bres}) is a minimal projective resolution of the trivial module ${}_B\k$.
\end{lem}

Next we compute the Yoneda product of $E(B)$. Note that $\Hom_B(B\ot C_{-n}\op B\ot C_{-n+1},{}_B\k)\cong C_{-n}^*\op C_{-n+1}^*$. For $\alpha\in C_{-n}^*$ and $\beta\in C_{-n+1}^*$, we view $(\alpha,\beta)$ as a homomorphism from $B\ot C_{-n}\op B\ot C_{-n+1}$ to ${}_B\k$.  Consider the following diagram:
{\small$$\xymatrix{\cdots B\ot C_{-n-k}\op B\ot C_{-n-k+1}\ar[r]^{\hspace{20mm}\delta^{-n-k}}\ar[d]_{g_k} &\cdots\ar[r]&B\ot C_{-n}\op B\ot C_{-n+1}\ar[d]_{g_o}\ar[dr]^{(\alpha,\beta)}\ar[r]&\cdots&\\
\cdots B\ot C_{-k}\op B\ot C_{-k+1}\ar[r]^{\hspace{20mm}\delta^{-k}} &\cdots\ar[r]& B\ar[r]^{\delta^{-1}}& \k\ar[r]&0,
}$$}
where the graded $B$-module homomorphisms $g_k$'s are defined as follows: for $k\ge1$,
$$    g_k(1\ot x_1\ot\cdots\ot x_k\ot x_{k+1}\ot\cdots\ot x_{n+k},0)
  =(1\ot x_1\ot\cdots\ot x_k\alpha(x_{k+1}\ot\cdots\ot x_{n+k}),0),
$$
and: {\small$$\begin{array}{l}
    g_k(0,1\ot x'_1\ot\cdots\ot x'_k\ot x'_{k+1}\ot\cdots\ot x'_{n+k-1})  \\
  =\left((-1)^k1\ot x'_1\ot\cdots\ot x'_k\beta(x'_{k+1}\ot\cdots\ot x'_{n+k-1}),1\ot x'_1\ot\cdots\ot x'_{k-1}\alpha(\sigma^{-1}(x'_{k})\ot\cdots\ot \sigma^{-1}(x'_{n+k-1}))\right);
  \end{array}$$}
for $k=0$,
$$ g_0(1\ot x_1\ot\cdots \ot x_{n},1\ot x'_1\ot\cdots\ot x'_{n-1})
=\alpha(x_1\ot\cdots \ot x_{n})1+\beta(x'_1\ot\cdots\ot x'_{n-1})1.
$$
A direct verification shows that the above diagram is commutative.

Since the projective resolution (\ref{Bres}) is minimal, we have $\Ext_B^n({}_A\k,{}_A\k)\cong C^*_{-n}\op C_{-n+1}^*$ for all $n\ge0$. Now assume $\alpha'\in C_{-k}^*$ and $\beta'\in C_{-k+1}^*$, we have:
\begin{equation}\label{eq1}
    (\alpha',\beta')\ast(\alpha,\beta)=(\alpha',\beta')\circ g_k=\left(\alpha'\cdot\alpha,\ (-1)^k\alpha'\cdot\beta+\beta'\cdot (\sigma^{-1})^!(\alpha)\right).
\end{equation}

\begin{prop}\label{extalg} Let $A$ be a Koszul algebra, $\sigma$ a graded automorphism of $A$, and $B=A[z;\sigma]$. Then $E(B)\cong \Gamma(A^!,{}_\epsilon A^!_\psi(-1))$, where $\psi=(\sigma^{-1})^!$ is the automorphism of $A^!$ dual to $\sigma^{-1}$.
\end{prop}
\proof Note that $A^!\cong C^\vee$ as graded algebras. The lemma is a direct consequence of the equation (\ref{eq1}). \qed

\section{Skew polynomial algebras with coefficients in Koszul Artin-Schelter regular algebras}\label{sec2}

In this section, $A$ is a Koszul Artin-Schelter regular algebra of global dimension $d$. The Artin-Schelter regularity of $A$ implies that $E(A)\cong A^!$ is a Frobenius algebra of length $d$ \cite{Smi}. Let $\varphi$ be the Nakayama automorphism of the Frobenius algebra $A^!$.

The following result was originally proved by Van den Bergh in \cite{VdB} in the Noetherian case. The result for general Koszul algebras was proved by Berger and Marconnet in \cite[Proof of Theorem 6.3]{BM}.

\begin{lem} \label{nakaut} Let $A=T(V)/(R)$ be a Koszul AS-regular algebra of global dimension $d$. Let $\varphi$ be the Nakayama automorphism of $A^!$, and $\phi:=\varphi^!$ the automorphism of $A$ dual to $\varphi$. Then $\Ext^i_{A^e}(A,A\ot A)=0$ for $i\neq d$, and $$\Ext_{A^e}^d(A,A\ot A)\cong A_\xi(d),$$ where $\xi=\epsilon^{d+1}\phi^{-1}$.
\end{lem}

The automorphism $\xi$ in the lemma is called the {\it Nakayama automorphism} of $A$.

The above lemma implies the following result (also see \cite{HVZ1}).

\begin{lem} \label{lem1} Let $A$ be a Koszul algebra. Then $A$ is Calabi-Yau if and only if $E(A)$ is a graded symmetric algebra.
\end{lem}

Now we may prove the main result of this section.

\begin{thm}\label{thm1} Let $A$ be a Koszul AS-regular algebra of global dimension $d$ with the Nakayama automorphism $\xi$. Then the skew polynomial algebra $B=A[z;\xi]$ is a Calabi-Yau algebra of dimension $d+1$.
\end{thm}
\proof Keep the notions as in Lemma \ref{nakaut}. By Proposition \ref{extalg}, we have $E(B)\cong\Gamma(A^!,{}_\epsilon A_{\psi}^!(-1))$, where $\psi=(\xi^{-1})^!$. Note that $\xi^{-1}=\epsilon^{d+1}\phi$. Then $\psi=(\xi^{-1})^!=\epsilon^{d+1}\varphi$. Therefore, we have ${}_\epsilon A_{\psi}^!\cong A^!_{\epsilon^d\varphi}$. Since $A^!$ is graded Frobenius with Nakayama automorphism $\varphi$, we have an isomorphism of $A^!$-bimodules $A^!_{\varphi}\cong (A^!)^*(-d)$, which implies $A^!_{\epsilon^d\varphi}\cong (A^!)^*_{\epsilon^d}(-d)$ since  $\epsilon\varphi=\varphi\epsilon$. Now we have:
\begin{equation}\label{eq5}
    E(B)\cong\Gamma(A^!,{}_\epsilon A_{\psi}^!(-1))\cong\Gamma(A^!,(A^!)^*_{\epsilon^d}(-d-1))=\Gamma(A^!,\epsilon^d,d+1).
\end{equation}
By Remark \ref{rem1}, $E(B)$ is a graded symmetric algebra. Lemma \ref{lem1} implies that $B$ is a Calabi-Yau algebra since $B$ is a Koszul algebra. \qed

\begin{cor}\label{corl} If $A$ is a Koszul Calabi-Yau algebra, so is $A[z]$.
\end{cor}

Let $V$ be a vector space of dimension $n$. Fix a basis $\{x_1,\dots,x_n\}$ of $V$, and $x^*_1,\dots,x^*_n$ the dual basis of $V^*$. Given an invertible $n\times n$ matrix $M$, let $f=(x_1,\dots,x_n)M(x_1,\dots,x_n)^t$, where the matrix multiplications are in the tensor algebra $T(V)$. The quadratic algebra $A$ has the following properties \cite{DV, Z}:

\begin{itemize}
\item[(i)] $A$ is a Koszul AS-regular algebra of global dimension 2;

\item[(ii)] $A$ is a domain;

\item[(iii)] the quadratic dual $A^!$ is defined by the matrix $M$ in the following way: there is a basis $\varpi$ of $A^!_2$, such that for $\alpha=a_1x^*_1+\cdots a_nx^*_n\in A^!_1$ and $\beta=b_1x^*_1+\cdots+b_nx^*_n\in A^!_1$, $\alpha\beta=\mathbf{a}M\mathbf{b}^t\varpi$, where $\mathbf{a}=(a_1,\dots,a_n)\in\k^n$ and $\mathbf{b}=(b_1,\dots,b_n)\in \k^n$ (cf. \cite{HVZ2});
\item[(iv)] the Nakayama automorphism of $A^!$ is defined in the way: $\varphi(\varpi)=\varpi$ and $$\varphi(x^*_1,\dots,x^*_n)=(x^*_1,\dots,x^*_n)M^{-1}M^t;$$
\item[(v)] the Nakayama automorphism $\xi$ of $A$ is defined in the following way: $$\xi(x_1,\dots,x_n)=-(x_1,\dots,x_n)M^tM^{-1}.$$
\end{itemize}

Let us check the generating relations of the skew polynomial algebra $B=A[z;\xi]$. Since $B$ is generated by $x_1,\dots,x_n,z$ and is quadratic, $B$ is defined by the relations $f=0$, $zx_1=\xi(x_1)z,\ \dots,\ zx_n=\xi(x_n)z$. Now $zx_i=-(x_1,\dots,x_n)M^tM^{-1}(0,\dots,1,\dots,0)^tz$. If we put $zx_1,\dots,zx_n$ into a column, we obtain
\begin{equation}\label{eq2}
    M^t\left(
    \begin{array}{c}
      zx_1 \\
      \vdots \\
      zx_n \\
    \end{array}
  \right)
=-M\left(
                                  \begin{array}{c}
                                    x_1z\\
                                    \vdots \\
                                    x_nz
                                  \end{array}
                                \right).
\end{equation}
Then we see that the skew algebra $B$ above is isomorphic to the algebra $B(f)$ constructed in \cite{BP}, which is the algebra generated by $x_1,\dots,x_n,z$ with relations $f=0$ and equations in (\ref{eq2}). As a corollary, we recover \cite[Theorem 2.10]{BP}.

\begin{thm}\cite{BP} The graded algebra $B(f)$ is a Calabi-Yau algebra of dimension 3.
\end{thm}

\begin{exa}{\rm Recall that a Noetherian AS-regular algebra $B$ of global dimension $d$ is called a {\it quantum polynomial algebra} if $B$ is a domain and has Hilbert series $H_B(t)=\frac{1}{(1-t)^d}$ (hence is Koszul).
Let $B$ be a Calabi-Yau quantum polynomial algebra of global dimension $d$. If $B$ is $\mathbb{Z}^2$-graded such that it is generated in degrees $(1,0)$ and $(0,1)$, and moreover $\dim B_{0,1}=1$, then $B=A[z;\xi]$, where $A=\op_{n\ge0}B_{n,0}$ is a quantum polynomial algebra of dimension $d-1$ \cite[Proposition 3.5]{KKZ}.
}\end{exa}

Let $A$ be a Koszul AS-regular algebra of global dimension 2. By \cite{DV, Z}, there is a finite dimension vector space $V$ with a fixed basis $\{x_1,\dots,x_n\}$ and an invertible $n\times n$ matrix $M$ such that $A\cong T(V)/(f)$ where $f=(x_1,\dots,x_n)M(x_1,\dots,x_n)^t$. We already know that the Nakayama automorphism of $A$ is defined by $\xi(x_1,\dots,x_n)=-(x_1,\dots,x_n)M^tM^{-1}$, and the Berger-Pichereau's algebra $B(f)\cong A[z;\xi]$. Let $M'$ be another invertible $n\times n$ matrix, and $f'=(x_1,\dots,x_n)M'(x_1,\dots,x_n)^t$. Let $A'=T(V)/(f')$. Denote by $\xi'$ the Nakayama automorphism of $A'$. The following result was proved in \cite[Theorem 3.4]{BP} (indeed, Berger-Pichereau did not assume that $M$ and $M'$ are invertible).

\begin{thm}\cite{BP} $B(f)\cong B(f')$ as graded algebras if and only if $M$ is congruent to a scalar multiple of $M'$; that is, there is an invertible $n\times n$ matrix $P$ and a scalar $k\in\k$ such that $M=kPM'P^t$. Moreover, if every element in $\k$ is a square in $\k$ then $B(f)\cong B(f')$ as graded algebras if and only if $M$ and $M'$ are congruent.
\end{thm}

However, we do not know whether there is a similar result for general Koszul AS-regular algebras.

\section{Superpotentials}

Let $V$ be a finite dimensional vector space. For the discussions in this section, we need additional notation. Let $\tau:V\ot V\longrightarrow V\ot V$ be the usual twisting map. For $d\ge2$, we set a sequence of maps: $\tau_d^0=1^{\ot d}:V^{\ot d}\to V^{\ot d}$, $\tau_d^1=\tau\ot 1^{\ot d-2}$, ... , $\tau^k_d=(1^{\ot k-1}\ot \tau\ot 1^{\ot d-k-1})\tau_d^{k-1}$ for all $k\ge2$.

Let $\sigma:V\to V$ be a linear bijective map. Recall that an element $w\in V^{\ot d}$ is called a {\it twisted superpotential of degree $d$} with respect to $\sigma$ if \begin{equation}\label{eq3}
    w=(-1)^{d-1}\tau_d^{d-1}\circ(\sigma\ot 1^{\ot d-1})(w).
\end{equation}
If $\sigma$ is the identity map, then $w$ is called a {\it superpotential}.

Let $\psi:V\to \k$ be a linear map, and let $u\in V^{\ot d}$. Following \cite{BSW}, we write: $$[\psi u]=(\psi\ot 1^{\ot d-1})(u)\text{, and }[u\psi]=(1^{\ot d-1}\ot\psi)(u).$$
More generally, if $\Psi\in (V^*)^{\ot k}$ $(k\leq d)$, we have: $$[\Psi u]=(\Psi\ot 1^{\ot d-k})(u),\text{ and } [u\Psi]=(1^{\ot d-k}\ot \Psi)(u).$$

One may check that an element $w\in V^{\ot d}$ is a twisted superpotential with respect to $\sigma$ if and only if, for all $\psi\in V^*$, $[\psi w]=(-1)^{d-1}[w(\psi\circ\sigma^{-1})]$.

For $\Psi\in {V^*}^{\ot k}$, define the partial derivation of a twisted superpotential $w$ to be $$\partial_\Psi(w)=[w\Psi].$$ Then $\partial_\Psi(w)\in V^{\ot d-k}$. The {\it derivation quotient algebra} $A(w,k)$ of $w$ is defined as follows \cite{BSW}: $$A(w,k)=T(V)/(\partial_\Psi(w):\Psi\in {V^*}^{\ot k}).$$
Since in this paper we only discuss the quadratic derivation quotient algebra, we simply write $A(w)$ for $A(w,d-2)$ for a twisted superpotential $w$ of degree $d$.

We now show that any twisted superpotential can be symmetrized into a superpotential by introducing an additional indeterminate. From the equation (\ref{eq3}), we have the following facts.

\begin{lem}\label{lem2} {\rm(i)} If $i\ge j\ge1$, then $\tau_d^i\circ\tau_d^j=\tau_{d}^{j-1}\circ(1\ot \tau_{d-1}^{i-1})$, and $\underbrace{\tau_d^{d-1}\circ\cdots\circ\tau_{d}^{d-1}}_{d\ factors}=1$;

 {\rm(ii)} Let $w\in V^{\ot d}$ be a twisted superpotential with respect to a bijection $\sigma$ of $V$. Then we have $$w=\sigma^{\ot d}(w).$$
\end{lem}
\proof (i) is trivial. For the statement (ii), we have
$$\begin{array}{ccl}
   w&=&(-1)^{d-1}\tau_d^{d-1}\circ(\sigma\ot 1^{\ot d-1})(w)\\
   &=& (-1)^{d-1}\tau_d^{d-1}\circ(\sigma\ot 1^{\ot d-1})((-1)^{d-1}\tau_d^{d-1}\circ(\sigma\ot 1^{\ot d-1})(w))\\
   &=&(-1)^{2(d-1)}\tau_d^{d-1}\circ\tau_d^{d-1}\circ(\sigma^{\ot 2}\ot 1^{d-1})(w)\\
   &\vdots&\\
   &=&(-1)^{d(d-1)}\underbrace{\tau_d^{d-1}\circ\cdots\circ\tau_d^{d-1}}_{d\ factors}\circ\sigma^{\ot d}(w)\\
   &=&\sigma^{\ot d}(w). \qed
 \end{array}
$$
\begin{prop}\label{prop} Assume that $w\in V^{\ot d}$ is a twisted superpotential with respect to a bijection $\sigma$ of $V$. We construct an element $\widehat{w}:=\widehat{w}(w,\sigma)\in (V\op \k z)^{\ot d+1}$ as follows: $$\widehat{w}:=\widehat{w}(w,\sigma)=\sum_{i=0}^d(-1)^i\tau_{d+1}^i(1\ot\sigma^{\ot i}\ot 1^{\ot d-i})(z\ot w).$$ Then $\widehat{w}$ is a superpotential of degree $d+1$.
\end{prop}
\proof We need to show the identity: $\widehat{w}=(-1)^d\tau_{d+1}^d(\widehat{w})$. This follows from the following computations:
$$\begin{array}{l}
\tau_{d+1}^d(\widehat{w})\\
= \displaystyle\tau_{d+1}^d\left(\sum_{i=0}^d(-1)^i\tau_{d+1}^i(1\ot\sigma^{\ot i}\ot 1^{\ot d-i})(z\ot w)\right)\\
= \displaystyle\tau_{d+1}^d(z\ot w)+\sum_{i=1}^d(-1)^i\tau_{d+1}^{i-1}(1\ot\tau_{d}^{d-1})(1\ot\sigma^{\ot i}\ot 1^{\ot d-i})(z\ot w)\\
= \displaystyle\tau_{d+1}^d(z\ot w)+\sum_{i=1}^d(-1)^i\tau_{d+1}^{i-1}(1\ot\tau_{d}^{d-1})(1\ot 1\ot\sigma^{\ot i-1}\ot 1^{\ot d-i})(1\ot\sigma\ot 1^{\ot d-1})(z\ot w)\\
= \displaystyle\tau_{d+1}^d(z\ot w)+\sum_{i=1}^d(-1)^i\tau_{d+1}^{i-1}(1\ot \sigma^{\ot i-1}\ot 1^{\ot d-i+1})(z\ot\tau_{d}^{d-1}(\sigma\ot 1^{\ot d-1})(w))\\
=\displaystyle\tau_{d+1}^d(z\ot w)+\sum_{i=1}^d(-1)^i\tau_{d+1}^{i-1}(1\ot \sigma^{\ot i-1}\ot 1^{\ot d-i+1})(z\ot(-1)^{d-1}w)\\
=\displaystyle\tau_{d+1}^d(z\ot w)+\sum_{i=1}^d(-1)^{d+i-1}\tau_{d+1}^{i-1}(1\ot \sigma^{\ot i-1}\ot 1^{\ot d-i+1})(z\ot w)\\
=\displaystyle\tau_{d+1}^d(z\ot w)+\sum_{j=0}^{d-1}(-1)^{d+j}\tau_{d+1}^{j}(1\ot \sigma^{\ot j}\ot 1^{\ot d-j})(z\ot w)\\
=\displaystyle\tau_{d+1}^d(z\ot \sigma^{\ot d}(w))+\sum_{j=0}^{d-1}(-1)^{d+j}\tau_{d+1}^{j}(1\ot \sigma^{\ot j}\ot 1^{\ot d-j})(z\ot w)\\
=(-1)^d\widehat{w}. \qed
  \end{array}
$$

Let $A=T(V)/(R)$ be a Koszul AS-regular algebra of global dimension $d$. It is established in \cite{BSW} and \cite{DV} that $A$ is defined by a twisted superpotential; that is, $A\cong A(w)$ for some twisted superpotential $w$ of degree $d$ with respect to a suitable bijection $\sigma:V\to V$. We give a ``visualized'' description of the bijection $\sigma$ and the twisted superpotential $w$.

Assume $\dim V=n$ and fix a basis $\{x_1,\dots,x_n\}$ of $V$. We use the definitions and notations as in Section \ref{sec}. Since $A$ is of global dimension $d$, we have $\dim C_{-d}=1$ and $\dim C_{-d+1}=n$. Choose a nonzero element $w\in C_{-d}$. Since $C_{-d}=\bigcap_{i=0}^{d-2}V^{\ot i} \ot R\ot V^{\ot d-i-2}$, it follows that $w\in V\ot C_{-d+1}\bigcap C_{-d+1}\ot V$. We fix a basis of $C_{-d+1}$, say, $\{\theta_1,\dots,\theta_n\}$. Since $w\in V\ot C_{-d+1}$, we may write $w$ as $w=(x_1,\dots,x_n)M(\theta_1,\dots,\theta_n)^t$ for some $n\times n$ matrix $M$ with entries in $\k$. On the other hand, $w\in C_{-d+1}\ot V$ implies that $w=(\theta_1,\dots,\theta_n)N(x_1,\dots,x_n)^t$ for some $n\times n$ matrix $N$. Let $\{x_1^*,\dots,x_n^*\}$ be the dual basis of $V^*$, and $\theta^*_1,\dots,\theta_n^*$ be the dual basis of $C_{-d+1}^*$. Since $A$ is AS-regular, $A^!\cong C^\vee$ is a graded Frobenius algebra of length $d$. For $\alpha=(x_1^*,\dots,x_n^*)(a_1,\dots,a_n)^t\in V^*=A^!_1$ and $\beta=(\theta_1^*,\dots,\theta_n^*)(b_1,\dots,b_n)^t\in C^*_{-d+1}=A^!_{d-1}$, it is easy to see that the Yoneda product is given by:
$$\alpha*\beta=(a_1,\dots,a_n)M\left(
                                 \begin{array}{c}
                                   b_1 \\
                                   \vdots\\
                                   b_n \\
                                 \end{array}
                               \right)w^*,
$$
and $$\beta*\alpha=(b_1,\dots,b_n)N\left(
                                 \begin{array}{c}
                                   a_1 \\
                                   \vdots\\
                                   a_n \\
                                 \end{array}
                               \right)w^*.$$
Now the Frobenius property of $A^!$ implies that both $M$ and $N$ are invertible matrices. Let $\varphi$ be the Nakayama automorphism of $A^!$. Then from the Yoneda product above, we see that $$\varphi(x_1^*,\dots,x_n^*)=(x_1^*,\dots,x_n^*)N^{-1}M^t.$$ Let $\phi:=\varphi^!$ be the automorphism of $A$ dual to $\varphi$. Then the Nakayama automorphism of $A$ is $\xi=\epsilon^{d+1}\phi^{-1}$, which acts on $A_1=V$ as follows: $$\xi(x_1,\dots,x_n)=(-1)^{d+1}(x_1,\dots,x_n)N^tM^{-1}.$$
We rewrite the element $w$ in terms of the Nakayama automorphism $\xi$ as follows:
$$\begin{array}{ccl}w&=&(\theta_1,\dots,\theta_n)N\left(
                                 \begin{array}{c}
                                   x_1 \\
                                   \vdots\\
                                   x_n \\
                                 \end{array}
                               \right)=(\theta_1,\dots,\theta_n)M^t(M^{-1})^tN\left(
                                 \begin{array}{c}
                                   x_1 \\
                                   \vdots\\
                                   x_n \\
                                 \end{array}
                               \right)\\
                               &=&(-1)^{d-1}(\theta_1,\dots,\theta_n)M^t\xi\left(
                                 \begin{array}{c}
                                   x_1 \\
                                   \vdots\\
                                   x_n \\
                                 \end{array}
                               \right).
                               \end{array}$$
That is, $w=(-1)^{d-1}\tau_d^{d-1}(\xi\ot 1^{\ot d-1})(w)$.

Summarizing the above arguments, we obtain the following lemma.

\begin{lem} Let $A=T(V)/(R)$ be a Koszul AS-regular algebra of dimension $d$ with $\xi$ the Nakayama automorphism. Then $w$ is a twisted superpotential with respect to $\xi|_V$.
\end{lem}

\begin{thm}\label{thm2} Let $A=T(V)/(R)$ be a Koszul AS-regular algebra of global dimension $d\ge2$ with $\xi$ the Nakayama automorphism. Then

{\rm(i)} $A\cong A(w)$, where $w$ is a nonzero element in $\bigcap_{i=1}^{d-2}V^{\ot i}\ot R\ot V^{d-2-i}$;

{\rm(ii)} $A[z;\xi]\cong A(\widehat{w})$, where the superpotential $\widehat{w}=\widehat{w}(w,\xi)$ is formed in Proposition \ref{prop}.
\end{thm}
\proof The statement (i) is essentially proved in \cite{DV} and \cite{BSW}. We include here our own proof. Assume $\dim R=m$. Since $A^!$ is Frobenius and $A^!_2=R^*$, we have that $A^!_{d-2}\cong C^*_{-d+2}$ is of dimension $m$. Fix a basis $\{r_1,\dots,r_m\}$ of $R$, and a basis $\{\vartheta_1,\dots,\vartheta_m\}$ of $C_{-d+2}$. As before, we let $\{r^*_1,\dots,r^*_m\}$ and $\{\vartheta^*_1,\dots,\vartheta^*_m\}$ be the dual bases of $R^*$ and $C_{-d+2}^*$ respectively. Note that we also have $w\in R\ot C_{-d+2}$. Hence there is an $m\times m$ matrix $L$ such that
\begin{equation}\label{eq4}
    w=(r_1,\dots,r_m)L(\vartheta_1,\dots,\vartheta_m)^t.
\end{equation}
For $\alpha=(r^*_1,\dots,r^*_m)(a_1,\dots,a_m)^t$ and $\beta=(\vartheta^*_1,\dots,\vartheta^*_m)(b_1,\dots,b_m)^t$, we have
$$\alpha*\beta=(a_1,\dots,a_m)L(b_1,\dots,b_m)^tw^*.$$ By the Frobenius property of $A^!$, we have that $L$ is invertible. Then from the expression of $w$ as in (\ref{eq4}), we see that $R=\{\partial_{\Psi}(w)|\Psi\in (V^*)^{\ot d-2}\}$. Therefore $A\cong A(w)$.

(ii) Since $w$ is a twisted superpotential with respect to $\xi$, $\widehat{w}$ is a superpotential of degree $d+1$ by Proposition \ref{prop}. Let $U=V\op\k z$. Then $\{x_0^*=z^*,x_1^*,\dots,x_n^*\}$ is a basis of $U^*$. Let us check the following facts:
\begin{itemize}
\item[(a)] $\{\partial_\Psi(\widehat{w})|\Psi\in (V^*)^{\ot d-1}\}=span\{z\ot x_i-\xi(x_i)\ot z|i=1,\dots,n\}$;

\item[(b)] $R=span \{\partial_\Psi(\widehat{w})|\Psi={x_{i_1}^*\ot\cdots\ot x^*_{i_{d-1}}} \text{ at least one of $i_1,\dots,i_{d-1}$ is zero} \}   $.
\end{itemize}
For $\Psi\in (V^*)^{\ot d-1}$, we have
$$\partial_\Psi(\widehat{w})=(1\ot 1\ot \Psi)[(z\ot w)-\tau_{d+1}^1\circ(1\ot\xi\ot 1^{\ot d-1})(z\ot w)].$$ Recall that $$w=(x_1,\dots,x_n)M\left(
                      \begin{array}{c}
                        \theta_1 \\
                        \vdots \\
                        \theta_n\\
                      \end{array}
                    \right).
$$ So, if we write $w=\sum_{i=1}^nx_i\ot y_i$, then $(y_1,\dots,y_n)=(\theta_1,\dots,\theta_n)M^t$. Since $M$ is invertible, we obtain that $y_1,\dots,y_n$ are linear independent in $V^{\ot d-1}$. Now we have $$\partial_\Psi(\widehat{w})=\sum_{i=1}^n (z\ot x_i-\xi(x_i)\ot z) \Psi(y_i).$$ Thus (a) follows.

For the identity (b), we choose $\Psi=x_{i_1}^*\ot\cdots\ot x_{i_{d-2}}^*\ot z^*$. Then $\partial_\Psi(\widehat{w})=(-1)^d(1\ot1\ot \Psi)(\xi^{\ot d}(w)\ot z)$. On the other hand, as we have seen that we may write $w$ in the form (\ref{eq4}). Again since $L$ is invertible, we have $$span\{\partial_\Psi(\widehat{w})|\Psi\text{ is of the form }{x_{i_1}^*\ot\cdots\ot x^*_{i_{d-2}}\ot z^*} \}=span\{(\xi\ot\xi)(r_i)|i=1,\dots,m\}.$$ As $A$ is Koszul and $\xi$ is the Nakayama automorphism of $A$, we have $span\{(\xi\ot\xi)(r_i)|i=1,\dots,m\}=R$. Since we obviously have $R\supseteq span\{\partial_\Psi(\widehat{w})|\Psi={x_{i_1}^*\ot\cdots\ot x^*_{i_{d-1}}} \text{ at least one of }\\ i_1,\dots,i_{d-1} \text{ is zero}\}$, (b) follows.

Finally, since $R+span\{z\ot x_i-\xi(x_i)\ot z|i=1,\dots,n\}$ is exactly the generating relations of $A[z;\xi]$, we have $A[z;\xi]\cong A(\widehat{w})$.\qed

\section{PBW-deformations}\label{pbw}

Let $A=A_0\op A_1\op\cdots$ be a positively graded algebra. Recall that a {\it PBW-deformation} of $A$ is a filtered algebra $U$ with an ascending filtration $0\subseteq F_0U\subseteq F_1U\subseteq F_2U\subseteq\cdots$ such that the associated graded algebra $gr(U)$ is isomorphic to $A$. If $A=T(V)/(R)$ is a Koszul algebra, then a PBW-deformation $U$ of $A$ is determined by two linear maps $\nu:R\to V$ and $\theta:R\to\k$ in the sense that $U\cong T(V)/(r-\nu(r)-\theta(r):r\in R)$ \cite{BG, PP}. If $\theta=0$, then we call $U$ an {\it augmented} PBW-deformation of $A$. The dual map $\nu^*$ of the linear map $\nu:R\to V$ induces a derivation $\delta_{A^!}$ on the dual algebra $A^!$ of $A$. If we view the linear map $\theta:R\to \k$ as an element in $A^!_2$, then $(A^!,\delta_{A^!}, \theta)$ is a cdga. We call $(A^!,\delta_{A^!},\theta)$ the {\it dual cdga} of $U$. Conversely, if there is a curved differential graded structure $(A^!,\delta_{A^!},\theta)$, then the dual map of the linear map $\delta_{A^!}|_{V^*}:V^*\to R^*$ and $\theta\in A^!_2=R^*$ define a PBW-deformation of $A$ \cite{PP}.

Now let $A$ be a Koszul AS-regular algebra of global dimension $d$, and let $U$ be a PBW-deformation of $A$. Assume that $\xi$ is the Nakayama automorphism of $A$. The following lemma was proved by Yekutieli \cite{Y} when $A$ is Noetherian. For the general case, see \cite{HVZ2}.

\begin{lem}\label{lem5} We have $\Ext_{U^e}^i(U,U\ot U)=0$ for $i\neq d$ and $\Ext_{U^e}^d(U,U\ot U)\cong U_\zeta$ as $U$-bimodules, where $\zeta$ is a filtration-preserving automorphism of $U$ such that $gr(\zeta)=\xi$.
\end{lem}

The automorphism $\zeta$ in Lemma \ref{lem5} is not unique. If $\zeta'$ is another automorphism of $U$ such that the conditions in the lemma hold, then $\zeta'$ differs from $\zeta$ by an inner automorphism, that is, there is a unit $u\in U$ such that for all $a\in U$, $\zeta'(a)=u\zeta(a)u^{-1}$. Hence $\zeta$ is unique up to inner automorphisms. We call an automorphism $\zeta$ satisfying the conditions in Lemma \ref{lem5} {\it a Nakayama automorphism} of $U$. Note that if $A$ is a domain, then there is a unique automorphism satisfies the condition in Lemma \ref{lem5}. Hence in this case, we may say ``the'' Nakayama automorphism of $U$.

Next we discuss the Calabi-Yau property of the skew polynomial algebra $U[z;\zeta]$ with $\zeta$ a Nakayama automorphism of $U$.

The skew polynomial algebra $U[z;\zeta]$ is also a filtered algebra with filtration: $F_0U[z;\zeta]=\k$, $F_nU[z;\zeta]=\sum _{i+j=n}F_iU z^j$ for all $n>0$ and $i,j\ge0$. It is easy to see that $gr(U[z;\zeta])\cong A[z;\xi]$. Hence we obtain:

\begin{lem} $U[z;\zeta]$ is a PBW-deformation of $A[z;\xi]$.
\end{lem}

The following result was proved in \cite{HVZ2}.

\begin{lem}\label{lem7} Let $B$ be a Koszul Calabi-Yau algebra of dimension $d$, and let $U=T(V)/(r-\nu(r)-\theta(r):r\in R)$ be a PBW-deformation of $B$, and let $(B^!,\delta_{B^!},\theta)$ be the cdga dual to $U$. If $\delta_{B^!}(B^!_{d-1})=0$, then $U$ is a Calabi-Yau algebra.

Conversely, if $B$ is a domain and $U$ is Calabi-Yau, then $\delta_{B^!}(B^!_{d-1})=0$.
\end{lem}

Let $B=A[z;\xi]$. As we have already proved in Section \ref{sec2}, $B$ is a Koszul Calabi-Yau algebra of dimension $d+1$. To see whether $U[z;\zeta]$ is a Calabi-Yau algebra, it is sufficient to see whether the condition $\delta_{B^!}(B^!_d)=0$ holds. By Proposition \ref{extalg}, $B^!\cong \Gamma(A^!,{}_\epsilon A^!_{\psi}(-1))$ where $\psi=(\xi^{-1})^!$. We need to write out the differential on $\Gamma(A^!,{}_\epsilon A^!_{\psi}(-1))$ induced by $\delta_{B^!}$ through the previous isomorphism. The following lemma is trivial.

\begin{lem}\label{lem6} Let $D$ and $D'$ be quadratic algebras. If there are invertible linear maps $f:D_1\to D'_1$ and $g:D_2\to D'_2$ such that the following diagram commutes: $$ \xymatrix{D_1\ot D_1\ar[d]^{\mu_D}\ar[r]^{f\ot f}&D'_1\ot D'_1\ar[d]^{\mu_{D'}}\\
D_2\ar[r]^g&D'_2,}$$ where the vertical maps are multiplications of $D$ and $D'$ respectively, then the linear map $f$ defines an isomorphism $\Phi:D\to D'$ in the following way: for any $x_1,\dots,x_n\in D_1$, $\Phi(x_1x_2\cdots x_n)=f(x_1)f(x_2)\cdots f(x_n)$.
\end{lem}

Let us write down $B^!_1$ and $B^!_2$ of $B^!$ explicitly. Write $\widehat{V}=V\op \k z$. As before, we fix a basis $\{x_1,\dots,x_n\}$ for $V$, and let $\{x_1^*,\dots,x_n^*\}$ be the dual basis of $V^*$. Let $\widetilde{r}_i=z\ot \xi^{-1}(x_i)-x_i\ot z$ for $i=1,\dots,n$, and $\widetilde{R}=span\{\widetilde{r}_1,\dots,\widetilde{r}_n\}\subseteq \widehat{V}\ot \widehat{V}$. Then $B=T(\widehat{V})/(\widehat{R})$, where $\widehat{R}=R\op \widetilde{R}$. Let $\{\widetilde{r}_1^*,\dots,\widetilde{r}_n^*\}$ be the dual basis of $\widetilde{R}$ and $z^*$ be the element in $\widehat{V}^*$ such that $z^*(z)=1$ and $z^*(V)=0$. We have $B^!_1=\widehat{V}^*$ and $B^!=\widehat{R}^*=R^*\op \widetilde{R}^*$, equivalently $B^!_1=A^!_1\op \k z^*$ and $B_2^*=A^!_2\op \widetilde{R}^*$.

Assume that the automorphism $\xi$ of $A$ acts on $A_1=V$ as:
\begin{equation}\label{eq11}
   \xi(x_1,\dots,x_n)=(x_1,\dots,x_n)P
\end{equation}
where $P=(p_{ij})$ is an invertible $n\times n$ matrix. Assume further $P^{-1}=(l_{ij})$. Then it is not hard to see that the product of two elements in $B^!_1$ is given as follows:  for $x_i^*, x_j^*\in A^!_1$, the product $x_i^*\cdot x_j^*$ is just the product in $A^!$; $z^*\cdot z^*=0$;

\centerline{$x_i^*\cdot z^*=-\widetilde{r}^*_i\in \widetilde{R}^*\subseteq B^!_2$ and $z^*\cdot x^*_i=\displaystyle\sum_{j=1}^n l_{ij}\widetilde{r}_j^*$.}

Recall that the graded algebra $B^!$ is isomorphic to $\Gamma(A^!,{}_\epsilon A^!_\psi(-1))$ with $\psi=(\xi^{-1})^!$ (Proposition \ref{extalg}). We construct an isomorphism from $\Gamma(A^!,{}_\epsilon A^!_\psi(-1))$ to $B^!$ in detail. Note that $\Gamma(A^!,{}_\epsilon A^!_\psi(-1))_1=A^!_1\op \k$ and $\Gamma(A^!,{}_\epsilon A^!_\psi(-1))_2=A^!_2\op V^*$. We define linear maps $f:\Gamma(A^!,{}_\epsilon A^!_\psi(-1))_1\to B^!_1$ and $g:\Gamma(A^!,{}_\epsilon A^!_\psi(-1))_2\to B^!_2$ as follows: $f(x^*_i,0)=x_i^*$ for all $i$ and $f(0,1)=z^*$; $g(\alpha,0)=\alpha$ for all $\alpha\in A^!_2$ and $g(0,x_i^*)=\widetilde{r}_i^*$ for all $i=1,\dots,n$. Now one may easily check that the conditions of Lemma \ref{lem6} above hold for $f$ and $g$. Therefore $f$ defines an isomorphism
\begin{equation}\label{eq9}
    \Phi:\Gamma(A^!,{}_\epsilon A^!_\psi(-1))\to B^!
\end{equation}
since both algebras are Koszul.

As before, let $\varphi$ be the Nakayama automorphism of $A^!$, and let $\phi=\varphi^!$ be the automorphism of $A$ dual to $\varphi$. Then $\xi=\epsilon^{d+1}\phi^{-1}$ by Proposition \ref{nakaut}. Hence $\varphi(x_i^*)=(-1)^{d+1}(\xi^{-1})^!(x_i^*)$ for all $i=1,\dots,n$. Since $A^!$ is Frobenius, there is an isomorphism of graded $A^!$-bimodules $\Theta:A^!_\varphi\to {A^!}^*(-d)$. Let $\varpi\in A^!_d$ be the element such that $\Theta(1)(\varpi)=1$. Then $\varpi$ is a basis of $A^!_d$. By the Frobenius property of $A^!$ again, we may choose elements $\omega_1,\dots,\omega_n$ in $A^!_{d-1}$ such that $x_i^*\omega_j=\delta^i_j\varpi$, where $\delta$ is the Kronecker delta function. Clearly, $\{\omega_1,\dots,\omega_n\}$ is a basis of $A^!_{d-1}$. Let $\varpi^*$ and $\{\omega_1^*,\dots,\omega^*_n\}$ be the dual basis of the space $(A^!_d)^*$ and $(A^!_{d-1})^*$ respectively. Consider the composition of the following isomorphisms:
$$h:{}_\epsilon A^!_\psi\overset{\epsilon^{d+1}}\longrightarrow {}_{\epsilon^d} A^!_\varphi\overset{\Theta}\longrightarrow {}_{\epsilon^d} {A^!}^*(-d).$$ We have $h(1)=\varpi^*$ and $h(x_i^*)=\sum_{j=1}^np_{ij}\omega^*_j$.
The isomorphism $h$ induces an isomorphism $\Gamma(A^!,{}_\epsilon A^!_\psi(-1))\longrightarrow\Gamma(A^!,{}_{\epsilon^d} {A^!}^*(-d-1))$. Combining this isomorphism with the inverse of $\Phi$ constructed in previous paragraph, we get an isomorphism of graded algebras:
\begin{equation}\label{eq10}
    \Psi:B^!\longrightarrow\Gamma(A^!,{}_{\epsilon^d} {A^!}^*(-d-1)).
\end{equation}
Now we have $\Psi(\alpha)=(\alpha,0)$ for all $\alpha\in A^!$, and
\begin{equation}\label{eq8}
    \Psi(z^*)=(0,\varpi^*)\text{ and } \Psi(\widetilde{r}_i^*)=(0,\sum_{j=1}^np_{ij}\omega^*_j)
\end{equation}
for all $i=1,\dots,n$.

Since $U[z;\zeta]$ is a PBW-deformation of $B=A[z;\xi]$, to study the curved differential structure of $B^!$, we need to pick a specific Nakayama  automorphism $\zeta$. The following result was proved in \cite{HVZ2}.

\begin{prop}\label{prop3} Let $A=T(V)/(R)$ be a Koszul AS-Gorenstein algebra of global dimension $d$, and let $A^!$ be its dual algebra. Assume that $\{x_1,\dots,x_n\}$ is a basis of $V$, and $\{x_1^*,\dots,x_n^*\}$ is the dual basis of $V^*$.

Let $U=T(V)/(r-\nu(r)-\theta(r):r\in R)$ be a PBW-deformation of $A$, and let $(A^!,\delta_{A^!},\theta)$ be the cdga dual to $U$. Choose a basis $\varpi$ of $A^!_d$, and assume that $\{\omega_1,\dots,\omega_n\}$ is the basis of $A^!_{d-1}$ such that $x_i^*\omega_j=\delta^i_j\varpi$. Assume further  $\delta_{A^!}(\omega_i)=\lambda_i\varpi$ for all $i=1,\dots,n$. Then $\Ext_{U^e}^i(U,U\ot U)=0$ for $i\neq d$, and  $$\Ext^d_{U^e}(U,U\ot U)\cong U_{\zeta},$$ where the automorphism $\zeta$ acts on the generator as follows:
\begin{equation}
\nonumber    \zeta(x_i)=\xi(x_i)+\lambda_i.
\end{equation}
\end{prop}

\noindent{\bf Convention:} From now on, $\zeta$ is the Nakayama automorphism of $U$ as defined in Proposition \ref{prop3}.

Note that in Proposition \ref{prop3}, we view $V\op \k$ as a subspace of $U$ through the obvious injective map. The scalars $\lambda_1,\dots,\lambda_n$ are independent of the choice of the basis $\varpi$. In fact, if we choose another element $\varpi'$ as a basis of $A^!_d$, then $\varpi'=k\varpi$ for some $k(\neq 0)\in \k$. Hence $x_i^*(k\omega_j)=\delta^i_j\varpi'$ for all $i,j=1,\dots,n$. Set $\omega'_i=k\omega_i$ for $i=1,\dots,n$. Then $\{\omega'_1,\dots,\omega'_n\}$ is the basis of $A^!_{d-1}$ satisfying the condition in the proposition. Clearly, we have $\delta_{A^!}(\omega'_i)=k\lambda_i\varpi=\lambda_i\varpi'$.

Now we can write down the linear maps that determine the PBW-deformation $U[z;\zeta]$ of $A[z;\xi]$. Recall that $A[z;\zeta]\cong T(\widehat{V})/(\widehat{R})$ with $\widehat{V}=V\op \k z$ and $\widehat{R}=R\op \widetilde{R}$.

\begin{lem}\label{lem8} $\widehat{U}:=U[z;\zeta]$ viewed as a PBW-deformation of $B=A[z;\xi]$ is determined by the following linear maps:
$$\begin{array}{cl}
 &\widehat{\nu}:R\op \widetilde{R}\to V\op \k z,\ \widehat{\nu}(r)=\nu(r)\ \text{for all $r\in R$}, and\ \widehat{\nu}(\widetilde{r}_i)=\lambda_i z\ (i=1,\dots,n);\\
&\widehat{\theta}:R\op \widetilde{R}\to \k,\ \widehat{\theta}(r)=\theta(r)\ \text{for all $r\in R$}, and\ \widehat{\theta}(\widetilde{r}_i)=0\ (i=1,\dots,n).
\end{array}$$ That is, $\widehat{U}\cong T(\widehat{V})/(\widehat{r}-\widehat{\nu}(\widehat{r})-\widehat{\theta}(\widehat{r}):\widehat{r}\in\widehat{R})$.
\end{lem}
\proof The lemma is clear since $z\xi^{-1}(x_i)-\lambda_i z=x_i z$ by Proposition \ref{prop3}, or equivalently, $z\xi^{-1}(x_i)-x_i z-\lambda_i z=0$ in $\widehat{U}$. \qed

Let us check the linear dual maps of $\widehat{\nu}$ and $\widehat{\theta}$. We have $\widehat{\nu}^*:\widehat{V}^*\longrightarrow \widehat{R}^*$, $\widehat{\nu}^*|_{V^*}=\nu^*$, and $\widehat{\nu}^*(z^*)=\sum_{i=1}^n\lambda_i \widetilde{r}^*_i$; $\widehat{\theta}^*=\theta^*:\k\to R^*\subseteq \widehat{R}^*$. Let $(A^!,\delta_{A^!},\theta_{A^!})$ be the cdga dual to the PBW-deformation $U$ of $A$. Since in the cdga $(B^!,\delta_{B^!}, \theta_{B^!})$, the differential $\delta_{B^!}$ is determined by $\widehat{\nu}^*$ and the curvature element $\theta_{B^!}=\widehat{\theta}$, we have $$\delta_{B^!}(\alpha)=\delta_{A^!}(\alpha)\text{ for all }\alpha\in V^*$$ and $$\delta_{B^!}(z^*)=\sum_{i=1}^n\lambda_i \widetilde{r}^*_i.$$

The cdga structure on $B^!$ induces a cdga structure on the graded algebra $\Gamma(A^!,{}_{\epsilon^d}(A^!)^*(-d-1))$ through the isomorphism $\Psi$ as in (\ref{eq10}) and (\ref{eq8}). Denote by $\Gamma:=\Gamma(A^!,{}_{\epsilon^d}(A^!)^*(-d-1))$. Let $(\Gamma,\delta_{\Gamma},\theta_{\Gamma})$ be the cdga induced by $(B^!,\delta_{B^!},\theta_{B^!})$. Then we have:
\begin{eqnarray}
\label{eq13}\theta_\Gamma&=&(\theta_{A^!},0);\\
\label{eq14}\delta_{\Gamma}(\alpha,0)&=&(\delta_{A^!}(\alpha),0), \text{ for all $\alpha\in A^!$}\\
 \label{eq15} \delta_{\Gamma}(0,\varpi^*)&=&\displaystyle\left(0,\ \sum_{j=1}^n\lambda_j(\sum_{i=1}^np_{ji}\omega^*_i)\right).
\end{eqnarray}

\begin{prop}\label{prop4}
Let $A=T(V)/(R)$ be a Koszul AS-regular algebra of global dimension $d$ with $\xi$ the Nakayama automorphism, and $A^!$ be its quadratic dual algebra with $\varphi$ the Nakayama automorphism. Let $U$ be a PBW-deformation of $A$, and $(A^!,\delta_{A^!},\theta_{A^!})$ be the cdga dual to $U$. If the composition $\epsilon^{d+1}\varphi$ is an automorphism of cdga $(A^!,\delta_{A^!},\theta_{A^!})$, then

{\rm(i)} The trivial extension $\Gamma_{cdg}(A^!, (A^!)^*(-d-1))$ of the cdga $(A^!,\delta_{A^!},\theta_{A^!})$ is the dual cdga of the PBW-deformation $U[z;\zeta]$ of $A[z;\xi]$, where $\zeta$ is the Nakayama automorphism of $U$ as in Proposition \ref{prop3};

{\rm(ii)} $U[z;\zeta]$ is a Calabi-Yau algebra. \end{prop}
\proof Keep the notions as before except we now choose the basis $\varpi$ in Proposition \ref{prop3} such that $\Theta(1)(\varpi)=1$. Let us check the differential $\delta^*_{A^!}$ on $(A^!)^*$. Recall that $\delta_{A^!}(\omega_i)=\lambda_i \varpi$ for all $i=1,\dots,n$ by assumption. Thus we have $$\delta_{A^!}^*(\varpi^*)=\sum_{i=1}^n\lambda_i\omega^*_i.$$
Note that $\varphi=\epsilon^{d+1}(\xi^{-1})^!$, and $\xi$ is represented as in (\ref{eq11}). An easy computation shows that $\varphi(\varpi)=\varpi$ and that:
$$\varphi(\omega_i)=(-1)^{d+1}\sum_{j=1}^np_{ji}\omega_j.$$
Now by assumption of the proposition $\delta_{A^!}\circ(\epsilon^{d+1}\varphi)=(\epsilon^{d+1}\varphi)\circ\delta_{A^!}$, we have $\delta_{A^!}^*\circ(\epsilon^{d+1}\varphi)^*=(\epsilon^{d+1}\varphi)^*\circ\delta_{A^!}^*$. Applying these morphisms to $\varpi^*$, we obtain:
$$(\epsilon^{d+1}\varphi)^*\circ\delta_{A^!}^*(\varpi^*)=(\epsilon^{d+1}\varphi)^*\left(\sum_{i=1}^n\lambda_i\omega^*_i\right)=\sum_{i=1}^n\lambda_i \sum_{j=1}^np_{ij}\omega^*_j;$$ and $$\delta_{A^!}^*\circ(\epsilon^{d+1}\varphi)^*(\varpi^*)=\delta_{A^!}^*(\varpi^*)=\sum_{i=1}^n\lambda_i\omega^*_i.$$ Hence we arrive at: \begin{equation}\label{eq12}
    \sum_{i=1}^n\lambda_i\omega^*_i=\sum_{i=1}^n\lambda_i \sum_{j=1}^np_{ij}\omega^*_j.
\end{equation}
Comparing the equations (\ref{eq15}) and (\ref{eq12}), we see that the differential $\delta_\Gamma$ on $\Gamma:=\Gamma(A^!,{}_{\epsilon^d}(A^!)^*(-d-1))$, induced from the cdga $(B^!,\delta_{B^!},\theta_{B^!})$, acts on the elements of degree 1 as:
\begin{eqnarray}
\label{eq16}\delta_{\Gamma}(\alpha,0)&=&(\delta_{A^!}(\alpha),0),\text{ for $\alpha\in A^!_1$} \\
\label{eq17}\delta_{\Gamma}(0,\varpi^*)&=&\displaystyle\left(0,\ \sum_{i=1}^n\lambda_i\omega^*_i\right).
 \end{eqnarray} Since $\Gamma$ is a quadratic algebra, the differential $\delta_\Gamma$ is determined by its action on the elements of degree 1.
By Lemma \ref{lem4}, the underlying graded algebra of the trivial extension $\Gamma_{cdg}(A^!, (A^!)^*(-d-1))$ of the cdga $(A^!,\delta_{A^!},\theta_{A^!})$ is exactly the graded algebra $\Gamma=\Gamma(A^!,{}_{\epsilon^d}(A^!)^*(-d-1))$. Let $\delta_{\Gamma_{cdg}}$ be the differential of $\Gamma_{cdg}(A^!, (A^!)^*(-d-1))$. By a straightforward check we have $\delta_{\Gamma_{cdg}}(\alpha,0)=(\delta_{A^!}(\alpha),0), $ for $\alpha\in A^!_1$ and $\delta_{\Gamma_{cdg}}(0,z^*)=\left(0,\ \sum_{i=1}^n\lambda_i\omega^*_i\right)$. Comparing these equations with (\ref{eq16}) and (\ref{eq17}), we see that the cdga $\Gamma_{cdg}(A^!, (A^!)^*(-d-1))$ is isomorphic to the cdga $(\Gamma,\delta_{\Gamma},\theta_{\Gamma})$. Hence the statement (i) follows.

Write $\widehat{\Gamma}:=\Gamma_{cdg}(A^!, (A^!)^*(-d-1))$. Then $\widehat{\Gamma}_{d}=A^!_d\op (A^!_1)^*$ and $\widehat{\Gamma}_{d+1}=\k$. Now it is clear that $\delta_{\Gamma_{cdg}}(\widehat{\Gamma}_{d})=0$. By Theorem \ref{thm1}, $A[z;\xi]$ is Calabi-Yau. Thus the statement (ii) follows from Lemma \ref{lem7}. \qed

As a special case of Proposition \ref{prop4}, we obtain the following theorem.

\begin{thm} \label{thm4} Let $A$ be a Koszul Calabi-Yau algebra of global dimension $d$, and let $U$ be an arbitrary PBW-deformation of $A$. Assume that $\zeta$ is the Nakayama automorphism of $U$ as in Proposition \ref{prop3}. Then $U[z;\zeta]$ is Calabi-Yau.
\end{thm}
\proof Since $A$ is Calabi-Yau, then the quadratic dual $A^!$ is graded symmetric; that is, the Nakayama automorphism of $A^!$ is $\varphi=\epsilon^{d+1}$. Then $\epsilon^{d+1}\varphi=id$, which is certainly an automorphism of the dual cdga $(A^!,\delta_{A^!},\theta_{A^!})$ of the PBW-deformation $U$. \qed

If $U$ is an augmented PBW-deformation, then the ground field $\k$ is a left $U$-module through the augmentation map. Let $E(U):=\op_{i\ge0}\Ext_U^i({}_U\k,{}_U\k)$. Note that the curvature element of the cdga $(A^!,\delta_{A^!},\theta_{A^!})$ dual to $U$ is zero. Thus $(A^!,\delta_{A^!})$ is a usual dga; that is, $\delta_{A^!}^2=0$. So, the cohomology $HA^!$ of $(A^!,\delta_{A^!})$ is a graded algebra.

\begin{prop} Let $A$ be a Koszul Calabi-Yau algebra, and $U$ an augmented PBW-deformation of $A$. Assume that $\zeta$ is the Nakayama automorphism of $U$ as in Propositon \ref{prop3}. Then $$E(U[z;\zeta])\cong\Gamma\left(H(A^!),{}_{\epsilon^{d+1}}H(A^!)^*(-d-1)\right).$$
\end{prop}
\proof If $U$ is an augmented PBW-deformation of $A$, then $U[z;\zeta]$ is an augmented PBW-deformation of $A[z]$. Hence $\Gamma_{cdg}(A^!, (A^!)^*(-d-1))$ is a dga. By \cite[Ch. 5, Proposition 6.1]{PP} and Proposition \ref{prop4}, $E(U[z;\zeta])$ is isomorphic to the cohomology algebra of the dga $\Gamma_{cdg}(A^!, (A^!)^*(-d-1))$. Now Lemma \ref{lem4} implies the desired isomorphism. \qed

In Proposition \ref{prop4}, we need the condition that the composition $\epsilon^{d+1}\varphi$ is an automorphism of the cdga $(A^!,\delta_{A^!},\theta_{A^!})$. Certainly, there is no reason to expect that $\epsilon^{d+1}\varphi$ is always compatible with the cdga structure on $A^!$. For example, if $A$ is an AS-regular algebra of global dimension 2, then $A^!$ is of length 2. Hence any linear map $\delta:A^!_1\to A^!_2$ and any element $\theta\in A^!_2$ form a cdga $(A^!,\delta,\theta)$. Below, we show that the condition that $\epsilon^{d+1}\varphi$ is compatible with the cdga structure on $A^!$ is necessary in case that $A$ is an AS-regular algebra of global dimension 2.

From now on, we assume that $A$ is an AS-regular algebra of global dimension 2. Then $A\cong T(V)/(f)$ where $V$ is an $n$-dimensional vector space with a fixed basis $\{x_1,\dots,x_n\}$, and $f=(x_1,\dots,x_n)M(x_1,\dots,x_n)^t$ with $M$ an invertible $n\times n$ matrix \cite{Z, DV}. Some properties of $A$ has been listed below Corollary \ref{corl}. In the following discussions, we keep the same notions as those in the items listed below Corollary \ref{corl}. We chose a new basis $\{\omega_1,\dots,\omega_n\}$ of $A^!_1$ as: $$(\omega_1,\dots,\omega_n)=(x_1^*,\dots,x_n^*)M^{-1}.$$ Then we have $x_i^*\omega_j=\delta_j^i\varpi$, where $\varpi$ is the basis of $A_2^!$ as in the item (iii) below Corollary \ref{corl}.

Let $U$ be a PBW-deformation of $A$. Since $A$ is a domain, as we pointed out earlier, there is a unique Nakayama automorphism. Hence we may say ``the'' Nakayama automorphism of $U$.

\begin{thm} \label{thm5} Let $U$ be a PBW-deformation of $A$ with $\zeta$ the Nakayama automorphism of $U$, and let $(A^!,\delta_{A^!},\theta_{A^!})$ be the dual cdga of $U$. Assume $\delta_{A^!}(\omega_1,\dots,\omega_n)=(\lambda_1\varpi,\dots,\lambda_n\varpi)$. The following are equivalent:
\begin{itemize}
\item[(i)] $\epsilon\varphi$ is an automorphism of the cdga $(A^!,\delta_{A^!},\theta_{A^!})$;
\item[(ii)] $U[z;\zeta]$ is Calabi-Yau;
\item[(iii)] $(\lambda_1,\dots,\lambda_n)M=-(\lambda_1,\dots,\lambda_n)M^t$.
\end{itemize}
\end{thm}
\proof That (i) implies (ii) follows from Proposition \ref{prop4}.

(ii)$\Longrightarrow$(iii). Let $B=A[z;\xi]$, $\widehat{U}=U[z;\zeta]$, and $(B^!,\delta_{B^!},\theta_{B^!})$ the dual cdga of $\widehat{U}$. Let $\Gamma:=\Gamma(A^!,(A^!)^*(-3))$, and $(\Gamma,\delta_{\Gamma},\theta_{\Gamma})$ the cdga induced by $(B^!,\delta_{B^!},\theta_{B^!})$ through the isomorphism $\Psi$ given in (\ref{eq10}) and (\ref{eq8}). The equation (\ref{eq15}) in the present case reads as follows:
\begin{equation}\label{eq18}
    \delta_{\Gamma}(0,\varpi^*)=(0,X),
\end{equation} where $X=-(\omega_1^*,\dots,\omega_n^*)(M^{-1})^tM(\lambda_1,\dots,\lambda_n)^t$. Now, since $\widehat{U}$ is Calabi-Yau and $B$ is obviously a domain, we have $\delta_{\Gamma}(\Gamma_2)=0$ by Lemma \ref{lem7}. Hence we have
$$\begin{array}{ccl}
0&=&\delta_{\Gamma}((\omega_i,0)*(0,\varpi^*))\\
&=&\delta_{\Gamma}(\omega_i,0)*(0,\varpi^*)-(\omega_i,0)*\delta_{\Gamma}(0,\varpi^*)\\
&=&(\lambda_i\varpi,0)*(0,\varpi^*)-(\omega_i,0)*(0,X)\\
&=&(0,\lambda_i\varpi\cdot \varpi^*)-(0,\omega_i\cdot X)\\
&=&(0,\lambda_i)-(0,\omega_i\cdot X),
 \end{array}
$$ where the notion ``$\cdot$'' in $\varpi\cdot \varpi^*$ and $\omega_i\cdot X$ is the left $A^!$-module action on $(A^!)^*$, and in the last identity, we identify $\k$ with $(A^!_0)^*$. Thus we obtain
\begin{equation}\label{eq19}
    \lambda_i=\omega_i\cdot X
\end{equation}
for all $i=1,\dots,n$. Note that $\omega_i\cdot X=-(0,\dots,0,1,0,\dots,0)(M^{-1})^tM(\lambda_1,\dots,\lambda_n)^t$. From Equation (\ref{eq19}), we obtain $(\lambda_1,\dots,\lambda_n)^t=-(M^{-1})^tM(\lambda_1,\dots,\lambda_n)^t$, and hence (iii) follows.

(iii)$\Longrightarrow$(i). We have
$$\begin{array}{ccl}
(\epsilon\varphi)\delta_{A^!}(\omega_1^*,\dots,\omega_n^*)&=&\epsilon\varphi(\lambda_1\varpi^*,\dots,\lambda_n\varpi^*)\\
&=&(\lambda_1\varpi^*,\dots,\lambda_n\varpi^*);
                                  \end{array}$$

and $$\begin{array}{ccl}
      \delta_{A^!}(\epsilon\varphi)(\omega_1^*,\dots,\omega_n^*)&=&-\delta_{A^!}\varphi(x_1^*,\dots,x_n^*)M^{-1}\\
      &=&-\delta_{A^!}(x_1^*,\dots,x_n^*)M^{-1}M^tM^{-1}\\
      &=&-(\lambda_1\varpi,\dots,\lambda_n\varpi)M^tM^{-1}.
      \end{array}
$$
Now the condition (iii) insures that $(\epsilon\varphi)\delta_{A^!}=\delta_{A^!}(\epsilon\varphi)$. Therefore $\epsilon\varphi$ is an automorphism of cdga $(A^!,\delta_{A^!},\theta_{A^!})$ since $A^!$ is of length 2. \qed

\vspace{5mm}

\subsection*{Acknowledgement}
The authors thank the referee for his/her valuable comments. The work is
supported by an FWO-grant and grants from NSFC (No. 11171067), ZJNSF (No. LY12A01013), Science and Technology Department of Zhejiang Province (No. 2011R10051), and SRF for ROCS, SEM.

\vspace{5mm}

\bibliography{}

\end{document}